\documentclass{article}
\usepackage{amsbsy,amssymb,amscd,amsfonts,latexsym,amstext,delarray,
amsmath}  \headsep=15pt
\def\pidr{ \pi_{1,DR}}
\def\piet{\pi_{1,\mbox{\'et}}}
\def\cP{ {\cal P} }
\def\Un{ \mbox{Un} }
\def\cQ{ {\cal Q} }
\def\cZ{ {\cal Z} }
\def\picr{ \pi_{1,cr} }
\def\cX{{\cal X} }
\def\cD{ {\cal D}}

\def\dim{{ \mbox{dim} }}
\def\Spec{{ \mbox{Spec} }}

\def\Hom{{ \mbox{Hom} }}
\def\AlgHom{ \mbox{Alg-Hom}}
\def\cA{ {\cal A} }

\def\Res{ \mbox{Res} }

\def\O{{ {\cal O} }}
\def\Aet{ A_{\mbox{\'et}} }
\def\ra{{ \rightarrow }}
\def\da{{ \downarrow }}

\def\a{{ \alpha }}

\def\g{{ \gamma }}
\def\d{{ \delta }}

\def\F{ {\bf F} }

\def\Un{ \mbox{Un} }
\def\Vect{ \mbox{Vect} }
\def\hra{{ \hookrightarrow }}
\def\da{{ \downarrow }}

\def\C{{ \mathbb{C} }}

\def\G{{ \Gamma }}
\def\Gal{{ \mbox{Gal} }}

\def\Z{{ \mathbb{Z}}}
\def\cC{ {\cal C} }
\def\Alb{ \mbox{Alb} }

\newtheorem{thm}{Theorem}

\newtheorem{lem}{Lemma}
\newtheorem{prop}{Proposition}

\usepackage{amsmath}
\usepackage{amsfonts}
\usepackage{amscd}
\usepackage{amssymb}
\def\Q{\mathbb{Q}}

\def\invlim{\varprojlim}

\def\cV{ {\cal V} }
\def\cD{ {\cal D} }
\def\P{ {\bf P}}
\def\An{ \mbox{An} }
\def\Col{ \mbox{Col} }
\def\Anloc{ \mbox{An}_{loc}}
\def\Anlog{ \mbox{An}_{\log}}
\def\Res{\mbox{Res}}

\title{The motivic fundamental group of ${\bf P}^1 \setminus \{0,1,\infty\}$ and the theorem of Siegel }
\author{ Minhyong Kim}

\begin{document}
\maketitle

One of the strong motivations for studying the arithmetic fundamental groups
of algebraic varieties comes from the hope that they will
provide group-theoretic access to Diophantine
geometry. This  is most clearly expressed by the so-called
`section conjecture' of Grothendieck: Let $F$ be a number field and
let $X \ra \Spec (F)$ be a smooth
hyperbolic curve. After some choice of base point, we get an
exact sequence
$$0\ra \hat{\pi}_1(\bar{X}) \ra \hat{\pi}_1(X)\stackrel{p}{ \ra }\G\ra 0$$
where $\G$ is the Galois group of $F$  and  $\hat{\pi}_1$ refers to the
profinite fundamental group. The section conjecture states that there
is a one-to-one correspondence between conjugacy classes of
splittings of $p$ and the geometric sections of the  morphism
$X \ra \Spec(F)$. There is apparently a compactness
argument that yields a direct implication from the section
conjecture to the  theorems of Siegel and Faltings.

The purpose of this paper is to illustrate  a somewhat different
methodology for deriving Diophantine consequences from a study of
the fundamental group.  To this end, we give a $\pi_1$ proof of
the theorem of Siegel on the finiteness of integral points for the
thrice-punctured projective line. The notion of a `$\pi_1$ proof'
may not be entirely well-defined, but it is hoped that the
techniques employed will make the meaning clear. Another way of
expressing the main idea is to say that we are using functions on
the fundamental group to prove Diophantine finiteness. Here, the
fundamental group refers to the unipotent motivic fundamental
group in the sense of Deligne \cite{deligne}. Although a rigorous
construction of such an object (for rational varieties) has  now
been given by Deligne and Goncharov  using Voevodsky's theory, our
proof does not require more than Deligne's original construction
using systems of realizations. More precisely, we will be using
the local and global \'etale fundamental groups $\piet$, the local
De Rham fundamental group $\pidr$, the crystalline fundamental
group $\picr$ and comparisons between them. Here and henceforward,
all $\pi_1$'s will denote $(\Q_p$-)unipotent completions, as they
are all we will be considering in this paper. In this connection,
it is rather striking that the motivic theory is capable of
yielding Diophantine finiteness, even though the motivic
fundamental group (at least the portion we use) could be viewed as
a cruder invariant of a variety than the pro-finite one. That is
to say,
 connections between the theory of motives  
and Diophantine geometry via $L$-functions is expected
in great generality. However, when viewed as invariants of
varieties, $L$-functions,
because they `factor through' the linear category of motives,
seem to provide
information in general
only about linearized invariants, e.g., Chow groups.
It is then natural that results about non-linear sets
 should evoke non-linear tools
like the fundamental group. What is somewhat surprising is
that even a mild degree of non-linearity (coming from the subcategory of
unipotent
group objects in the category of motives) 
can still provide substantial information.

Here is a slightly more precise outline of the proof: In the
discussion above, set $F=\Q$, and
$X=\P^1_{\Q }\setminus \{0,1,\infty\}$. Let $S$ be a finite set of primes.
Fix any $p\notin S$ and an $S$-integral point
$x$ of $\cX:=\P^1_{\Z}\setminus \{0,1,\infty\}$. Let $T=S\cup \{p\}$.
Denote by  $Y$ the reduction of $\cX$ mod $p$. Let $y\in Y$ be the reduction
mod $p$ of the point $x$. We will use a $p$-adic {\em unipotent
Albanese map}
$$UAlb_x:  X(\Q_p) \cap ]Y[\ra \pidr (X,x)(\Q_p)$$
associated to the basepoint $x$ and defined on the tube of $Y$
inside $X(\Q_p)$, that is, the points that reduce mod $p$ to
points of $Y$. This map is constructed by considering the class of
the compatible pair of torsors of paths $\pidr(X\otimes \Q_p
,x,x')$ and $ \picr(Y,y,y')$ associated to a point $x'$ with
reduction $y'\in Y$. These are torsors for the crystalline and De
Rham fundamental groups $\pidr(X\otimes \Q_p,x)$ and $
\picr(Y,y)$, respectively. A simple classification of such
compatible pairs of torsors allows us to canonically associate to
the pair a point in $\pidr (X\otimes\Q_p,x)$ which we define to be
$UAlb_x (x')$. Now $\pidr (X\otimes \Q_p,x)$ is a pro-unipotent
algebraic group over $\Q_p$ with a coordinate ring $A_{DR}$ that
is generated as a $\Q_p$-vector space by functions $\a_w$ indexed
by the words $w$ in two letters $A,B$. An important point is that
the functions $g_w:=\a_w\circ UAlb$ are restrictions to
$X(\Q_p)\cap ]Y[$ of {\em Coleman functions} on $X(\C_p)$ and
hence, have nice analytic properties. We will show that they are
$\C_p$-linearly independent, so no non-zero function from
$\cA_{DR}$ pulls back to the zero function on $X(\Q_p)$. On the
other hand, when we examine the image of the integral points
$\cX(\Z_S)$ under the Albanese map, we find that it is essentially
contained inside the image of another `map' $$C:H^1_f(\G_T,
\piet(X,x))\ra \pidr (X\otimes \Q_p,x)$$ from a suitable
continuous global cohomology set to the De Rham fundamental group.
This `map' is algebraic and is obtained from global-to-local
restriction and $p$-adic Hodge theory. We explain the quotation
marks: In fact, $H^1_f(\G_T, \piet(X,x))$ has the natural structure
of a pro-algebraic variety.  If we look at various quotients
$[\pidr]_n$ and $[\piet]_n$ with respect to the descending central
series of these fundamental groups, what we actually have are
finite-dimensional varieties and algebraic maps
$$C_n:H^1(\G_T, [\piet(X,x)]_n)\ra [\pidr (X\otimes \Q_p,x)]_n$$
whenever $p$ is sufficiently large with respect to $n$ (we will
explain this in detail in section 3).  We can also consider the
level n unipotent Albanese maps
$$UAlb_n:X(\Q_p)\cap ]Y[ \ra [\pidr(X,x)]_n(\Q_p)$$ obtained
from $UAlb$ via composition with the natural projections, and we
find that
$$UAlb_n(\cX (\Z_S)) \subset \mbox{Im}(C_n).$$
 However, the descending central
series filtration on $H_f^1(\G_T, [\piet(X,x)]_n)$ together with a
Galois cohomology computation of Soul\'e allows us to get explicit
bounds on the dimensions of these Galois cohomology varieties. The
upshot then is that for large $n$ and $p$, the image of $H_f^1(\G_T,
[\piet(X,x)]_n)$ under the map to $ [\pidr (X\otimes \Q_p,x)]_n$
lies inside a proper subvariety. Therefore, some non-zero element
of $A_{DR}$ vanishes on this image, and hence, on the image of
$\cX (\Z_S)$. This fact, together with
 the identity principle for
Coleman functions and
compactness yields the finiteness of Siegel's theorem.

The reader familiar with the method of  Chabauty
(\cite{ch},\cite{coleman}) will immediately recognize our proof to be a
`non-abelian lift' of his. As such, we believe that many
generalizations and refinements should be possible and hope to
discuss them in the near future. The present paper, however, was
motivated by the wish to work out in full detail one non-trivial
example, thereby testing the strength of the techniques involved
in carrying out this lift.

\section{Torsor spaces}

We need some elementary preliminaries on topological vector
spaces. Let $B$ be a complete Hausdorff topological field of
characteristic zero. Given a finite-dimensional vector space $V$
over $B$, there is a unique topology on $V$  compatible with the
vector space structure. This can be described by choosing any
isomorphism $V\simeq B^n$ and using the product topology. We will
be considering $B$-vector spaces $R$, possibly
infinite-dimensional. It will be convenient to topologize such $R$
by giving them the inductive limit topology coming from the family
of all finite-dimensional subspaces. Thus, a map $f:R\ra A$ is
continuous if and only if
 its restriction to any finite-dimensional subspace $V\subset R$
is continuous. The inclusion $V\hra R$ of a finite-dimensional
subspace is then automatically continuous. In fact, it is a direct
consequence of the definitions that any $B$-linear map is
continuous: Say $f:R_1 \ra R_2$ is $B$-linear. Take $V_1\subset
R_1$ finite-dimensional. Then $V_2:=f(V_1)\subset R_2$ is also
finite-dimensional. Since $f|V_1$ can be factored as $V_1\ra V_2
\hra R_2$ and both arrows are continuous, we are done. This
argument is typical of those involving the inductive limit
topology. Also obvious is that any vector subspace is closed. Note
that the topology is Hausdorff: Let $v,w \in R, v\neq w$ and let
$V\subset R$ be the subspace generated by $v, w$. By choosing a
basis of $R$, we can construct a projection $p:R\ra V$ which must
be continuous. Now find $O_v, O_w \subset V$, disjoint open
subsets of $V$ containing $v$ and $w$ respectively. Then
$p^{-1}(O_v)$ and $p^{-1}(O_w)$ separate $v$ and $w$.

The inductive limit topology can be applied, in particular, to the
situation where $R$ is a $B$-algebra, or to any affine $n$-space
$R^n$ over $R$. Suppose $V\subset R^n$ is finite-dimensional. Then
each of the projections $V_i:=p_i(V)$ to the components is
finite-dimensional and $V\subset \prod_i V_i$. So
finite-dimensional subspaces of this product form are co-final.
Hence, it suffices to check continuity of maps on such subspaces.
In fact, it is clearly sufficient to consider subspaces of the
form $W\times W \times \cdots \times W$ with $W$
finite-dimensional in $R$.

\begin{lem}
Any polynomial map $f:R^n \ra R^k$ is continuous.
\end{lem}
{\em Proof.} Consider a subspace $\prod_iW \subset R^n$ as above.
Let $d$ be the maximal degree of the monomials occurring in any
component of $f$ and let $\{a_j\}$ be the set of coefficients of
$f$. Let $\{w_1,\ldots, w_m\}$ be a basis for $W$ and consider the
$B-$subspace $V$ of $R$ generated by all the
$a_jw_1^{\a_1}w_2^{\a_2}\cdots w_m^{\a_m}$ as $\alpha=(\a_1, \a_2,
\ldots, \a_m)$ runs over all multi-indices of weight $\leq d$.
Obviously, $f$ takes $\prod W$ to $\prod V$ and is continuous on
$\prod W$. Since the inclusion $\prod V \hra R^n$ is continuous,
we are done. $\Box$

Given any affine $R$-scheme $X$ of finite-type, we can give the
$R$-points $X(R)$ of $X$ the topology induced by any embedding
$X\hra A^n_R$. The lemma above shows that this topology is
independent of the embedding. Also, we get
\begin{lem}
If $X\ra Y$ is a map of affine $R$ schemes of finite-type, then
the induced map on $R$ points in continuous.
\end{lem}

Equally obvious from the definitions is the
\begin{lem} If $R\ra T$ is a map of $B$-algebras,
then $X(R)\ra X(T)$ is continuous for any affine $R$-scheme $X$
of finite-type.
\end{lem}

For the cohomological considerations below, it will be useful to
note the following
\begin{lem}
Let $R$ be a $B$-vector space with the inductive limit topology
and let $C\subset R$ be compact. Then $C$ is contained in a
finite-dimensional subspace $V\subset R$.
\end{lem}
{\em Proof.} We cannot have $C=R$, since then, $C\cap Bv=Bv\simeq
B$ for some one-dimensional $Bv \subset R$ would be compact, and
there is no infinite compact-field (recall that $B$ has
characteristic zero). By choosing $v\notin C$ and translating to
$C-v$, we can assume that $0\notin C$. Suppose we had an infinite
collection $\{v_n\}_{n=1}^{\infty} \ldots$ of linearly independent vectors in
$C$. By passing to a subsequence, can assume that the $v_n$
converge to $v\in C$. Write $R=W_1\oplus W_2$ where $W_1$ is the
subspace generated by the $v_i$'s and $W_2$ is a complement. We
can write  $v=\Sigma_ic_i v_i+w$, where $w\in W_2$. Only finitely
many $v_i$ occur in the sum, say $1\leq i \leq N$. Let $p_2$ be
the projection to $W_2$. Then $0=p_2(v_n)$ converges to $p_2(v)=w$, so
$w=0$. Now write $W_1=W_1'\oplus W_1''$
where $W_1'$ is generated by $v_1,\ldots, v_N$ and
$W_1''$ is generated by $v_{N+1},\ldots$.
Consider the projection $p:W_1\ra W_1'$ determined by this decomposition.
 Then $p(v_n)$ converges to
$p(v)=\Sigma_ic_i v_i$. But $p(v_i)=0$ for $i>N$. So $\Sigma_ic_i
v_i=0$. Therefore, we conclude that $0=v\in C$, a contradiction.
$\Box$

We wish to consider certain continuous cohomology spaces for a
compact topological  group $G$. The setting will be slightly more
general than appears necessary for this paper since we wish
also to look ahead to future work.
 Assume that $G$ acts continuously
on $B$ by field automorphisms and denote by $F$ the fixed field
for this action. In particular, the action could be trivial and
$B=F$. We also denote by $K$ a subfield of $F$ such that $[F:K]$
is finite. Thus, $F$ and $K$ are both complete with the induced
topology. For any $K$-algebra $R$, we get a $B$-algebra
$B\otimes_KR$ which we equip with a $G$-action via $g(b\otimes
r)=(gb)\otimes r$. This action is clearly continuous.

For each $K$-algebra $R$, we have  a category $\cC_R$ whose
 objects are
modules $\cP$ over $B\otimes_K R$ equipped with the following
data:

(1) A continuous, semi-linear action of $G$. Here, the continuity
means that the action map $G\times \cP \ra \cP$ is continuous,
while the semi-linearity refers to $g(bx)=g(b)g(x)$ for $b\in
B\otimes R, x\in \cP$. Note that $\cP$ can simply be regarded as a
$B$-vector space, and hence, the topology is that discussed above.

(2) An increasing filtration $W$ indexed by $\Z$ such that

-$W_n\cP=0$ for $n< 0$ and $\cup_nW_n\cP=\cP$.

-each $W_n\cP$ is stable under the $G$ action and is  a finitely
generated module over $B\otimes_KR$.

 The
morphisms in this category consist of

$$\Hom(\cP,\cQ)=\invlim_n \Hom (W_n\cP,W_n\cQ)$$

where $\Hom (W_n\cP,W_n\cQ)$ is the set of linear maps of
$B\otimes_KR$-modules that commute with the $G$-action.

 There is a natural tensor product in
this category obtained by putting the tensor-product filtration on
$\cP\otimes_{B\otimes_K R} \cQ$:
$$W_n(\cP\otimes \cQ)=\Sigma_{i+j=n}W_i\cP\otimes W_j\cQ$$

Let  $\cP$ be a finitely generated $B$ algebra in the category
$\cC:=\cC_K$. Therefore, the structure maps are all required to be
$\cC$-morphisms. If we denote $P=\Spec(\cP)$, we have a natural
action of $G$ on the set of points $P(B\otimes R)$. Explicitly, an
element $g\in G$
 takes an $B-$algebra homomorphism
$\phi: \cP \ra B\otimes R$ to $g\phi g^{-1}$.

\begin{lem} This $G$-action is continuous.
\end{lem}
Take $n$ large enough so that $W_n\cP$ has a set of algebra
generators for $\cP$.
 Then for any $B$-algebra $S$, we have
$$P(S) \hra \Hom_{B}(W_n\cP, S)=\Hom_{B}(W_n\cP,B)\otimes_{B} S$$
But $\Hom_{B}(W_n\cP,B)$ is just a finite-dimensional vector space
over $B$, so this gives us an embedding
 $P\hra \Hom_{B}(W_n\cP,B)$ into a finite-dimensional
affine space over $B$
 that is compatible with the
$G$-action on $B\otimes_KR$-points. That is, $$P(B\otimes_K R)
\subset \Hom_{B}(W_n\cP,B)\otimes_B B\otimes _KR =
\Hom_{B}(W_n\cP,B)\otimes _KR $$ with the induced topology.
Therefore, we need only check the continuity of the action on
$\Hom_{B}(W_n\cP,B)\otimes_KR $. For any finite dimensional
$K$-subspace $V$ of $R$, $\Hom_{B}(W_n\cP,B)\otimes_K V$ is a
$G$-invariant subspace of $\Hom_{B}(W_n\cP,B)\otimes_KR$ and these
subspaces are cofinal among finite-dimensional $B$-subspaces of
$\Hom_{B}(W_n\cP,B)\otimes_K R$. The map $G\times
(\Hom_{B}(W_n\cP,B)\otimes_K V)\ra \Hom_{B}(W_n\cP,B)\otimes_K V$
is continuous by the continuity of the original $G$-action on
$W_n\cP$ and on $B$. Therefore, $G\times
(\Hom_{B}(W_n\cP,B)\otimes_K R) \ra \Hom_{B}(W_n\cP,B)\otimes_K R$
is also continuous. $\Box$

We will be interested in a pro-unipotent algebraic group $U$ in
the dual category $\cC^{\circ}$. Thus, in $\cC$, the corresponding
object is a $B$-algebra $\cA$ with the structure of a Hopf
algebra.
 $\cA$ is equipped with a multiplication
$$m:\cA\otimes \cA \ra \cA,$$
a comultiplication
$$\d: \cA \ra \cA\otimes \cA,$$
a unit $1 \in \cA$, a counit $e:\cA \ra B$ and an antipode $i:\cA
\ra \cA$ which are morphisms in the category $\cC$ and are
compatible in the usual sense. Let $\cA_1\subset \cA_2 \subset
\ldots \subset \cA$ denote  the filtration of $\cA$ by subalgebras
corresponding to the descending central series of $U:=\Spec(\cA)$,
normalized so that $U_1=U/[U,U]=\Spec(\cA_1)$. Let
$U_n=\Spec(\cA_n)$. We assume that each $U_n$ is a unipotent
algebraic group over $B$, i.e., that each $\cA_i$ is a finitely
generated integral domain. By a $U$-torsor, we will mean an affine
$\cC$-scheme $P$ (in the sense of Deligne (\cite{deligne}, 5.4
))with an action of $U$ that makes it into a $U$-torsor in the
usual sense (loc. cit.). However, we recall that the structure
maps of the action are required to be in the category $\cC$.
Therefore, the coordinate ring $\cP$ of $P$ is an object of $\cC$
and we are given a map $a: \cP \ra \cP \otimes \cA$ in the
category $\cC$ that induces a free and transitive group action on
points.  In our definition, we also assume that the torsors are
strictly compatible with the weight filtration, in the sense that
$\dim W_n\cP=\dim W_n\cA$ for each $n$. Now, note that when we
forget the Galois action, such a torsor is always trivial, i.e.,
has a $B-$rational point $x\in P$, since the group $U$ is
unipotent. Such a point gives rise to an isomorphism of Hopf
algebras:
$$\cP \stackrel{a}{\ra }
\cP\otimes \cA \stackrel{ x \otimes 1}{\ra} \cA$$ The isomorphism
preserves the filtration $W$:
$$W_n(\cP) \ra \Sigma_{i+j=n}W_i(\cP)\otimes W_j(\cA) \ra \Sigma_{j\leq n}
W_j(\cA)\subset W_n(\cA)$$ since the filtration is increasing. In
fact, this map is {\em strict} for the filtration, i.e., induces
an isomorphism $W_n\cP \ra W_n\cA$ for each $n$ by the equality of
dimensions.

The torsor $P$ gives rise to a
collection of $U_n$-torsors $P_n$ obtained by push-out:
$$P_n:=(P\times U_n)/U$$
where $U$ acts on the product via the diagonal action, that is,
$g(p,x)=(pg,g^{-1}x)$.

More generally, we can consider a $U_R:=\Spec (\cA
\otimes_KR)$-torsor $P=\Spec (\cP)$ in the category $\cC_R^o$,
defined in an obvious way analogous to the above discussion. The
only further requirement is that $W_n\cP$ is a free
$B\otimes_KR$-module of finite rank equal to $\dim_B W_n \cA=
\mbox{rank}_{B\otimes_KR} \cA \otimes_K R$. In
this case as well, the choice of a point $x\in P(B\otimes R)$
determines an isomorphism
$$\cP \simeq \cA \otimes_KR$$
that preserves the weight filtration. By changing base to closed
points $y$ of $\Spec (B\otimes R)$ whereupon we get torsors $P_y$
for $(U_R)y$ over the fields $k(y)$, we see that $W_n\cP\otimes
k(y)\simeq W_n \cA \otimes_K R \otimes_{B\otimes R} k(y)=W_n \cA \otimes_B k(y)$ by
dimension considerations, and hence, that $W_n\cP \simeq W_n
\cA\otimes_KR$.

The basic classification goes as follows:

\begin{prop}
The isomorphism classes of $U_R$ torsors are in bijection with the
continuous cohomology set
$$H^1(G, U(B\otimes_KR)).$$
\end{prop}

The continuous cohomology occurring in the statement is defined in
the standard way  (\cite{serre}, VII Appendix)
 which we  review briefly. Given any $K$-algebra
$R$ we extend scalars to the $B$-algebra $B\otimes_K R$ and  give
$U(B\otimes_K R)$ the topology and $G$-action discussed above.
This gives rise to the set $C^i(G,U(B\otimes_K R))$ of  continuous
i-cochains, which are defined to be continuous maps $c:G^i \ra
U(B\otimes_K R)$. The boundary maps $d:C^i(G,U(B\otimes_K R)) \ra
C^{i+1}(G,U(B\otimes_K R))$ are defined in a standard way at least
for $i=0$ and $i=1$. We recall the explicit description. In degree
0, $C^0(G,U(B\otimes_K R))=U(B\otimes_K R)$ and for $u\in
U(B\otimes_K R)$, we have $$(du)(g)=ug(u^{-1}).$$ For $c:G\ra
U(B\otimes_K R)$ a continuous map,
$$dc(g_1,g_2)=c(g_1g_2)(g_1c(g_2))^{-1}c(g_1)^{-1}.$$
All the $C^i$ are pointed sets, where the point is the constant
map taking values in the identity element $e$ of $U(B\otimes_K
R)$. $e$ will also be used to denote any of the corresponding
cochains. We then have the continuous 1-cocycles
$Z^1(G,U(B\otimes_K R))\subset C^1(G,U(B\otimes_K R))$ defined as
$d^{-1}(e)$.
 Thus, it consists of continuous
maps $c:G \ra U(B\otimes_K R)$ such that
$c(g_1g_2)=c(g_1)g_1c(g_2)$. Consider then the action of
 $U(B\otimes_K R)$  on $Z^1(G,U(B\otimes_K R))$ by
$(uc)(g)=uc(g)g(u^{-1})$. We define
$$H^1(G, U (B\otimes_K R)):=U(B\otimes_K R))\backslash Z^1(G,U(B\otimes_K R))$$
In the case where $V$ is a vector group over $B$, we also have
conventional definitions of $Z^i(G,V(B\otimes_K R))$ (the
$i-$cocycles), $B^i(G,V(B\otimes_K R))$ (the $i-$coboundaries),
and $$H^i(G,V(B\otimes_K R))=B^i(G,V(B\otimes_K R))\backslash Z^i(G,V(B\otimes_K
R))$$ for all $i$ (\cite{serre}, VII.2).

{\em Proof of Proposition.} Let $P=\Spec (\cP)$ be a $U_R$-torsor.
Thus, $\cP$ is an object of the category $\cC_R$ and we are given
a $\cC$-morphism
$$a:\cP \ra \cP \otimes \cA$$
specifying the $U$-action. Now choose a point $x\in P(B\otimes
R)$. Given an element $g\in G$, we have a unique element $u_g \in
U(B\otimes R)$ such that $gx=xu_g$. We get thereby a map $G\ra
U(B\otimes R)$. To check continuity of this map we give an
alternative description. $x$ is an algebra homomorphism $\cP \ra
B\otimes R$. Hence, as described above, we can form the composite
$$s_x:\cP \stackrel{a}{\ra} \cP \otimes_{B\otimes R} (\cA \otimes_K R)\stackrel{x\otimes 1}{\ra} \cA\otimes_KR$$
which is a continuous isomorphism of $G\otimes_KR$-algebras. For
each $g\in G$, we then have a continuous isomorphism of algebras
$c_x(g):= s_x g s_x^{-1}g^{-1}: \cA\otimes_KR  \ra \cA\otimes_KR$.
As above, let $n$ be large enough so that $W_n\cA$ contains a
generating set for $\cA$. Then $c_x(g)$ is determined by its
restriction to $W_n\cA\otimes_KR$. Since $W_n\cA\otimes_KR$ and
$W_n\cP$
 are finite-rank over
$B\otimes_KR$ and the original $G$-action on either side is
continuous, the map $g \mapsto c_x(g)|W_n\cA\otimes R$ is clearly
continuous. On the other hand, the elements $y\in U(B\otimes_KR)$
act on $W_n\cA\otimes R$ (and $\cA\otimes R$) by
$$\phi_y: W_n\cA\otimes R \stackrel{m}{\ra} W_n (\cA \otimes_KR \otimes_{B\otimes_KR} \cA\otimes_KR) \stackrel{y^{-1}\otimes 1}{\ra}
W_n\cA\otimes_KR$$ and this determines an affine embedding of $U$
from which $U(B\otimes_KR)$ gets the induced topology.

Claim: $\phi_{u_g}=c_x(g)|W_n\cA\otimes_KR$

Let $f\in \cA\otimes_KR$. We compute $c_x(g)(f)$ on points (an
obvious dual argument gives the desired equality): For a point
$y\in U$ and $h\in \cA\otimes R$, the $G-$actions are related by
$(gh)(y)=g(h(g^{-1}y))$. Also, given $h\in \cA \otimes_KR$,
$s_x^{-1}(h)(z)=h(y)$ where $y$ is determined by $z=xy$.
Then the proof of the claim is an exercise in careful bracketing:
$$\begin{array}{rl}
c_x(g)(f)(y)=(s_x(g(s_x^{-1}(g^{-1}f))))(y)=&
(g(s_x^{-1}( g^{-1}f)))(xy)\\
=g((s_x^{-1}( g^{-1}f))(g^{-1}(x)g^{-1}(y)))&=g((g^{-1}f)(w)) \ \ (\mbox{for}\ xw=g^{-1}(x)g^{-1}(y))\\
=g(g^{-1}(f(g(w))))&=f(g(w))
\end{array}$$
But $xy=g(x)g(w)=xu_gg(w)$ so $g(w)=u_g^{-1}y$. That is to say,
the end result is $(\phi_{u_g}f)(y)$, as desired.

Now, if $g_1,g_2 \in G$, then
$$g_1g_2x=g_1(xu_{g_2})=g_1(x)g_1(u_{g_2})=xu_{g_1}g_1(u_{g_2})$$
so that $g\mapsto u_g$ is a 1-cocycle. It is straightforward to
check at this point that different choices of $x$ give us
equivalent cocycles, and hence, we get a well-defined class $c(P)
\in H^1(G, U(B\otimes_K R))$.

In the other direction, given a 1-cocycle $c \in Z^1(G,
U(B\otimes_KR))$, we can use it to twist the Galois action on
$\cA\otimes_KR$ by letting $\rho_c(g)(x)= \phi_{c(g)}gx$. The
usual formula
$$\begin{array}{rl}
\phi_{c(g_1g_2)}g_1g_2(x)=\phi_{c(g_1)g_1c(g_2)}g_1g_2(x)=& \\
\phi_{c(g_1)}\phi_{g_1(c(g_2))}g_1g_2(x))=\phi_{c(g_1)}g_1\phi_{(c(g_2))}g_1g_2(x))&
=\phi_{c(g_1)}g_1(\phi_{(c(g_2))}g_2(x))\end{array}$$ shows that
this is a group action. The continuity of the cocyle gives the
continuity of the action. Since the action by $\phi_{c(g)}$
preserves the filtration $W$, so does the representation $\rho_c$.
Denote by $\cP(c)$ the filtered algebra $\cA\otimes_KR$ with this
twisted action. We give $P(c)=\Spec(\cP(c))$ the $U_R$-action by
using the group law on $U_R$:
$$m:\cP (c) \ra \cP(c) \otimes_{B\otimes_KR} \cA\otimes_KR$$
We need to check that this is compatible with the Galois action.
Note that $\phi_y=(y^{-1}\otimes 1)\circ m$ for any $y\in
U(B\otimes_KR)$. That is to say, we need to check the identity
$$[(c(g)\otimes 1 \circ m)\otimes 1](g\otimes g) (m(a))=
m( (c(g)\otimes 1 \circ m)(ga)$$ for all $g \in G$ and $a \in
\cA\otimes_KR$. But $(g\otimes g)\circ m=m\circ g$ since $m$ is
$G$-equivariant for the original action. Thus, we have to check
$$[(c(g)\otimes 1 \circ m)\otimes 1]\circ m (a)) =m\circ
(c(g)\otimes 1 \circ m)(a)$$ for all $a\in \cA$. We can again
check this by dualizing an argument on points of $U\times U$: The
left-hand side evaluated on $(h,k)$ becomes $f((gh)k)$ while the
right-hand side is $f(g(hk))$.

The usual computation (cf. \cite{serre}, X.2) shows that equivalent cocycles
give isomorphic actions and that the two correspondences are
inverses to each other. $\Box$

If $R\ra S$ is a map of $K$-algebras, we have the induced map
$$H^1(G, U(B\otimes_K R)) \ra H^1(G, U(B\otimes_K S))$$
which we view therefore as defining a functor on $K$-algebras.
That is, we define the functor $H^1(G,U)$ by
$$H^1(G, U)(R):=H^1(G, U(B\otimes_K R))$$
for any $K$-algebra $R$. We can also define similar functors
$C^1(G,U)$, $Z^1(G, U)$ and for vector groups $V$, $C^2(G,V)$,
$Z^2(G,V)$ and $B^2(G, V)$. We denote $U^1=U$ and
$U^{i+1}=[U,U^i]$. Also, $U_i:=U/U^{i+1}$.

We make the following important assumption:

\medskip

$H^j(G, U^i/U^{i+1}(B))$ are finite-dimensional $F$-vector
spaces for each $i, j$.

\medskip

In section three we will also encounter the natural condition
that the inclusion $\cA^G\hra \cA$ of the $G$-invariants induces an
isomorphism $(\cA^G \otimes_F B)\simeq \cA$ compatible with the
$G$-action.

In this case, we
denote by $U^G$, the unipotent algebraic group $\Spec(\cA^G)$ and we
see that for any $K$-algebra $R$, we have
$$H^0(G,U)(R)=\AlgHom_G[\cA, B\otimes_KR],$$ where
$\AlgHom$ refers to $B$-algebra homomorphisms while the subscript
$G$ restricts to the $G$-invariant ones. However,
$$\begin{array}{rl}
\AlgHom_G[\cA, B\otimes_K R]=&
\AlgHom_G[(\cA)^G\otimes_FB, B\otimes_KR]
\\
=\AlgHom_G[(\cA)^G, B\otimes_K R]=&
\AlgHom[(\cA)^G, (B\otimes_KR)^G]=\\
=\AlgHom[(\cA)^G,  F\otimes_K R]=& U^G(F\otimes_K R)
\end{array}$$
That is, the functor $H^0(G,U)$ is represented by the Weil
restriction $\Res_{F/K}(U^G)$.
 However, this last condition is not necesary for the following:
\begin{prop} Suppose $H^0(G, U^i/U^{i+1})(K)=0$ for each $i$.
Then the functor $H^1(G,U)$ is representable by an affine pro-algebraic
variety over $K$.
\end{prop}
{\em Proof.}
Note that the assumption implies that
$H^0(G,U^i/U^{i+1})(R)=0$ for all $i$ , and hence that
$H^0(G,U_i)(R)=0$, given any $k$-algebra $R$.

First we prove the elementary fact
that  each of the $H^1(G,U^i/U^{i+1})$ are
representable. $V:=U^i/U^{i+1}$ is a finite-dimensional vector
group and we have $V(B\otimes_KR)=V(B)\otimes_KR$. We claim that
the natural map
$$f: C^n(G,V(B))\otimes_KR \ra C^n(G, V(B)\otimes R)
\ra C^n(G, V(B\otimes R))$$ is an isomorphism. To check
injectivity, let $\{r_i\}$ be a $K$-basis for $R$ and let
$c=\Sigma c_i \otimes r_i \in C^n(G,V(B))\otimes_KR $. Suppose
$f(c)=0$. Then $f(c)(g)=\Sigma c_i(g)\otimes r_i=0$ for all $g\in
G$. Since $r_i$ form a basis, this implies that $c_i(g)=0$ for
each $i$ and $g$, so $c_i=0$, and hence, $c=0$. To check
surjectivity, let $c:G^n \ra V(B)\otimes R$ be continuous. Since
$G^n$ is compact, the image has to lie inside a finite-dimensional
$B$-subspace. If the subspace is generated by $\Sigma_i
b_{ij}\otimes r_i$, then at most finitely many $r_i$'s occur in
all these sums, so the image of $c$ lies inside a subspace of the
form $V(B)\otimes_K W$ where $W\subset R$ has finite
$K$-dimension. Let $\{w_i\}$ be a basis for $W$. Then $c(g)=\Sigma
c_i(g) \otimes w_i$ defines a finite set of $c_i:G^n \ra V(B)$ and
the element $\Sigma c_i\otimes w_i \in C^n(G,V(B))\otimes W$ such
that $f(\Sigma c_i\otimes w_i )=c$. By the exact same argument, we
can check that $Z^1(G,V(B\otimes R))=Z^1(G,V(B))\otimes R$.

Now since $H^1(G,V)$ is defined by the exact sequence
$$V(B\otimes R) \ra Z^1(G,V(B\otimes R)) \ra H^1(G,V(B\otimes R)) $$
we get
$$H^1(G,V)(R)\simeq H^1(G,V(B))\otimes_K R$$
in a way functorial in $R$. So $H^1(G,V)$ is represented by the
finite-dimensional vector group $H^1(G,V(B))$.

Now we will prove the theorem  inductively for $H^1(G,U_n)$.

We have an exact sequence of algebraic groups
$$0\ra U^{n+1}/U^{n+2} \ra U_{n+1} \ra U_{n} \ra 0$$
which realizes $U_{n+1}$ as a $U^{n+1}/U^{n+2}$-torsor over $U_n$
(and as a central extension of algebraic groups). But
$U^{n+1}/U^{n+2}$ is a vector group and $U_n$ is affine so this torsor
splits. Choose
 an algebraic splitting $s:U_n \ra U_{n+1}$
for the projection $U_{n+1} \ra U_n$. This induces a continuous
map $U_{n}(S) \ra U_{n+1}(S)$ for any $B$-algebra $S$ which is in
fact functorial in $S$. Thus, for any $K$ algebra $R$, we have a
split  exact sequence
$$0\ra U^{n+1}/U^{n+2}(B\otimes R) \ra U_{n+1}(B\otimes R) \ra U_n(B\otimes
R) \ra 0$$

Let $n=1$. Then $U_1=U/[U,U]$ is a vector group. We have shown
that $H^1(G,U_1)$ is representable by a vector group.

Assume we have proved the representability for $n$. Consider the
surjective map
$$Z^1 (G, U_n(B\otimes R)) \ra H^1(G,U_n)(R)$$
Taking for $R$ the coordinate ring of $H^1(G,U_n)$, we have the
element of $H^1(G,U_n)(R)$ corresponding to the identity map.
Choosing a lifting to $Z^1 (G, U_n(B\otimes R))$ gives us a
functorial splitting
$$i:H^1(G,U_n)\ra Z^1 (G, U_n)$$
Composing with the section $s$ and the boundary map $d:C^1(G,U_{n+1})
\ra C^2(G,U_{n+1})$, we get a map $dsi: H^1(G,U_n) \ra Z^2(G,
U^{n+1}/U^{n+2})$. Define $I(G,U_n):=(dsi)^{-1}(B^2(G,U^{n+1}/U^{n+2}))$.
Note that $dsi$ composed with the natural quotient map $Z^2(G,
U^{n+1}/U^{n+2}) \ra H^2(G, U^{n+1}/U^{n+2})$ realizes the connecting
homomorphism  $$\d:H^1(G,U_n)\ra H^2(G, U^{n+1}/U^{n+2})$$ in a functorial way (cf. \cite{serre}, VII, Appendix,
Prop. 2) and $I=\d^{-1}(0)$, so it
is a closed affine subvariety of $H^1(G,U_n) $. The proof given
above shows that
$$C^1(G, U^{n+1}/U^{n+2}(B\otimes R))\simeq
C^1(G, U^{n+1}/U^{n+2}(B))\otimes R$$
and
$$ B^2(G, U^{n+1}/U^{n+2}(B\otimes R))\simeq B^2(G, U^{n+1}/U^{n+2}(B))\otimes R$$
so if we choose a $K$-linear splitting
$$a:B^2(G, U^{n+1}/U^{n+2}(B))\ra C^1(G, U^{n+1}/U^{n+2 }(B))$$
of the boundary map, then we get a functorial splitting
$$
a:B^2(G, U^{n+1}/U^{n+2})\ra C^1(G, U^{n+1}/U^{n+2})$$ Thus, we can
define
$$b(x)=(si)(x)(adsi)(x)^{-1}$$
to get a map $b:I(G,U_n) \ra Z^1(G,U_{n+1})$. By composing with
the quotient map $Z^1(G,U_{n+1})\ra H^1(G,U_{n+1})$ we get a
functorial section
$$I(G,U_n) \ra H^1(G,U_{n+1})$$
of the surjection
$$H^1(G,U_{n+1}) \ra I(G,U_n)$$

{\bf Claim:} For each $R$, we have an exact sequence
$$0\ra H^1 (G,U^{n+1}/U^{n+2})(R)\ra
H^1(G,U_{n+1})(R) \ra I(G,U_n)(R)\ra 0$$
in the sense that the left hand group acts
freely on the middle set, and the surjection identifies
$I(G,U_n)(R)$ with the set of orbits.

{\em Proof of Claim:}
The only non-evident part is the freeness of the action.
So assume we have $c_{n+1} \in Z^1(G,U_{n+1})(R) $
that maps to $c_n \in Z^1(G,U_{n})(R)$. Suppose
$v\in Z^1 (G,U^{n+1}/U^{n+2})(R)$ stabilizes the class
of $c_{n+1}$. Then there exists a $u\in U_{n+1}$
such that $c_{n+1}(g)v(g)=u^{-1} c_{n+1}(g) g(u)$ for
all $g\in G$. Projecting to $U_n$ gives
$c_n(g)=\bar{u}^{-1}c_n(g)g(\bar{u})$ for each $g$,
where $\bar{u}$ is the projection of $u$ to $U_n$.
Hence, $c_n(g)g(\bar{u})c_n(g)^{-1}=\bar{u}$.
By induction on $n$ and our assumption that $H^0(G,U^i/U^{i+1})=0$
for all $i$, this implies that $\bar{u}=0$. Therefore,
$u\in U^{n+1}/U^{n+2}$ and hence, $v(g)=u^{-1}g(u)$
is a coboundary. $\Box$

The claim together with the section constructed
above induces an isomorphism of functors
$$H^1(G, U^{n+1}/U^{n+2})\times I(G,U_n) \simeq H^1(G,U_{n+1})$$
and concludes the proof that each $H^1(G,U_n)$ is represented by
an affine variety. Now, the surjectivity of the map
$U_{n+1}(B\otimes_KR)\ra U_n(B\otimes_KR)$ can be used to show
easily that $H^1(G,U)=\invlim H^1(G,U_n)$ as set-valued functors.
$\Box$

We will be considering in section 3 the important situation
where we compare cohomology sets over $K$ and $B$. That is
to say, suppose $U$ is defined over $K$ and assume that the
finite-dimensionality
assumption preceding Prop.2 is satisfied over both $K$ and $B$.
 Given any $K$-algebra $R$, $U_n(R)$
acts as a subgroup on $U_n(B\otimes_KR)$ and gives rise to
 the `exact sequence'
$$0\ra U_n(R) \ra U_n(B\otimes_KR) \ra 
U_n(B\otimes_KR)/U_n(R)\ra 0$$
from which we get a small part of a long exact sequence of pointed sets
$$H^0(G, 
U_n(B\otimes_KR)/U_n(R))\ra H^1(G, U_n(R)) \ra H^1(G,U_n(B\otimes_KR))$$ 
The topology on the quotient set will not be too important for
us since we will be considering only $H^0$'s for general unipotent
groups. For example, it is entirely straightforward to check
that one does get the portion of the long exact sequence displayed above
with the $H^1$'s being continuous cohomology.
However, we will need to consider the topology somewhat in the case
of a vector group $V$ over $K$. The inclusion
$K \hra B$ gives $K$ the induced topology and makes $B$
into a topological $K$-vector space. Since $K\subset B$ is
of course a finite-dimensional $K$-subspace, it has a topological
complement $C\subset B$, which is a closed $K$-subspace
such that $B=K\oplus C$ (\cite{koethe}, section 10.7 (8)).
Now we give $V(B)/V(K)$ the quotient topology and $V(B\otimes R)/V(R)
=(V(B)/V(K))\otimes R$ the inductive limit of the topology
coming from the subspaces $(V(B)/V(K))\otimes W$, where $W\subset R$
is finite-dimensional. Then we see that the exact sequence
$$0\ra V(R) \ra V(B\otimes R) \ra V(B\otimes R)/V(R) \ra 0$$
has a continuous $K$-linear splitting, giving rise to a long exact
sequence of continuous cohomology groups. In particular,
we see that $H^1(G,V(B)/V(K))$ is finite-dimensional
and that 
$$H^1(G,V(B\otimes R)/V(R))\simeq H^1(G,V(B)/V(K))\otimes R$$
i.e., the functor
$$H^1(G, V^B/V)(R):=H^1(G,V(B\otimes R)/V(R))$$
is represented by the vector group $ H^1(G,V(B)/V(K))$.

\begin{prop} The functor $$H^0(G,  U_n^B/U_n):
R\mapsto H^0(G, U_n(B\otimes_KR)/U_n(R))$$
is represented by an affine variety over $K$.
\end{prop}

{\em Proof.}
As before, if we examine the case of a vector group $V$,
we get an exact sequence
$$H^0(G, V(B)) \ra H^0(G, V(B)/V(K)) \ra H^1(G, V(K))$$
that exhibits $H^0(G,V(B)/V(K))$ as a finite-dimensional
$K$-vector space that represents the functor in question.
Of course, in the vector group case, the functor
$H^1(G, V^B/V)$ is also
representable by the vector space $H^1(G,V(B)/V(K))$
as discussed above.

The general case is again proved by induction on $n$.
That is, setting $V:=U^{n+1}/U^{n+2}$,
 we easily verify that we have an exact sequence of
$G$-sets
$$0\ra  V(B\otimes R)/V(R) \ra 
 U_{n+1}(B\otimes R)/U_{n+1}(R) \ra
U_{n}(B\otimes R)/U_{n}(R)\ra 0$$
from which we get a functorial exact sequence
$$\begin{array}{l}
0\ra H^0(G, V(B\otimes R)/V(R) \ra H^0(G,
U_{n+1}(B\otimes R)/U_{n+1}(R))
\ra H^0(G, U_{n}(B\otimes R)/U_{n}(R)) \\
\stackrel{\delta}{\ra} H^1(G, 
 V(B\otimes R)/V(R))
\end{array}
$$
Taking the inverse image of the basepoint under the connecting
map again defines a subfunctor
$$I_n(R)=\delta^{-1} (0) \subset H^0(G,U_{n}(B\otimes R)/U_{n}(R))$$
When we apply the inductive hypothesis to regard
$H^0(G,  U_n^B/U_n)$ as an affine variety, then $I_n$
becomes a closed subvariety. Denoting by $A(I_n)$
the coordinate ring of $I_n$, we then have an element in
$I_n(A(I_n))$ corresponding to the identity, which can then be
lifted to $H^0(G,U_{n+1}^B/U_{n+1})(A(I_n))$.
This determines a splitting $I_n\ra H^0(G,U_{n+1}^B/U_{n+1})$
and hence,    an isomorphism of functors
$$H^0(G,  U_{n+1}^B/U_{n+1}) \simeq H^0(G,  U_n^B/U_n)\times I_n$$
that gives us the desired representability. $\Box$

Clearly the map
$$H^0(G, U_n(B\otimes_KR)/U_n(R))\ra H^1(G, U_n(R))$$
is functorial and hence is given by an algebraic
map of varieties (when the latter is also representable as in Prop.2). 
Thus, taking $R=K$,
we see that the inverse image of the basepoint under the
map $H^1(G,U(K))\ra H^1(G,U(B))$ is the set of $K$-points
of a subvariety $H_f^1(G,U_n)$ of $H^1(G,U_n)$.

\section{The p-adic K-Z equation}

In this section, $k$ denotes an algebraic closure of $\F_p$.
$\O_F$ is a complete discrete valuation ring with residue field
$k$, and $F$ denotes the fraction field of $\O_F$.
 $\Vect_F$
refers to the category of finite-dimensional vector spaces over
$F$.
 $\cX$ denotes the projective line minus the
three points $0,1,\infty$: $\cX=\Spec \O_F[z][1/z(z-1)]$. $X$ is
the generic fiber of $\cX$ and $Y$ the special fiber. Following
Deligne, we define $\pidr(X,x)$, the De Rham fundamental group of
$X$ with basepoint $x\in X(F)$ to be the Tannakian fundamental
group of the category $\Un$ of unipotent vector bundles with
connection on $X$ associated to the fiber functor $e_x: \Un \ra
\Vect_F$ that takes $(V,\nabla)$ to the vector space $V_x$
(\cite{deligne}, , 10.27).  This definition can be generalized  to include
tangential basepoints as well as to  torsors of paths
$\pidr(X,x,y)$ associated to two points $x,y \in X$
(\cite{deligne}, 15.28). That is, if let $x\in \P^1\setminus X$ and let $v$
be an $F$-rational tangent vector to $\P^1$ at $x$ we get a fiber
functor $e_v$ according to the following procedure. First we
canonically extend $(V,\nabla)$ to a connection
$(\bar{V},\bar{\nabla})$ on $\P^1$ with log poles along
$\{0,1,\infty\}$ and nilpotent residues.
 Deligne then describes a procedure for associating to this data
 a connection $\Res(V,\nabla)$ on $T_x(\P^1)\setminus \{0\}$.
 $\Res(V,\nabla)$ is functorial in $(V,\nabla)$ and  taking
 the fiber of $\Res(V,\nabla)$ at $v$ defines the functor
 $e_v$.
 Now given two of these (possibly tangential)
 fiber functors $e_x$ and $e_y$, we get the functor of
 isomorphisms from $e_x$ to $e_y$, which is a right torsor
 $\pidr (X,x,y)$ for the group $\pidr (X,x)$.
Define $U:=\pidr(X,v)$ where $v$ is the tangent vector $d/dz$ at
$0 \in {\bf P}^1$. Define $P(x):=\pidr(X,v,x)$ for $x\in X(F)$. As
in the previous section, we will also consider the lower central
series $U_n$ and the associated torsors $P_n(x)$.

Let $X_{an}$ denote the rigid analytic space associated to $X$.
 We will be using the ring $\Col_a$
 of Coleman functions on $X_{an}$
with respect to the embedding $j:Y\hra \P^1_k$ (\cite{besser}, section 4).
For any rigid analytic space $Z$, we will denote by $\An (Z)$
the ring of rigid analytic functions on $Z$.
 Given  $a\in F$, denote by $\log_a$ the branch
 of the p-adic log normalized by the condition
$\log_a(p)=a$. Given a point $y\in \P^1(\bar{k})$, we have the
tube $]y[\subset \P^1_F$ consisting of points that reduce to $y$.
For $0\leq r<1$,
 $X_{an}(r)$ denotes the rigid analytic space obtained
by removing from
 $\P^1_{an}$ all closed disks of radius $r$ around
 $0,1,$ and $\infty$.
Define $\Anloc^a:=\Pi_{y\in \P^1(\bar{k})} \Anlog^a(y)$ where
$\Anlog^a(y):=\An(]y[)$ if $y\neq 0,1, \infty$ while
$\Anlog^a(y):=\lim_{r\ra 1}\An(]y[\cap X_{an}(r))[\log_a z_y]$ if
$y=0,1,$ or $\infty$ \cite{furusho}.  Here, $z_y$ denotes $z$,
$z-1$, or $1/z$ for $y=0,1,$ or $\infty$, respectively. Denote by
$\An^{\dagger}$ the ring $\Gamma (\P^1_{an},
j^{\dagger}\O_{\P^1_{an}})$, where $j^{\dagger}$ is Berthelot's
dagger functor \cite{berthelot}. The ring $\Col_a$ naturally
contains $\An^{\dagger}$ and is equipped with a map $r_y:\Col_a\ra
\Anlog(y)$ for each $y$. We refer to Besser (\cite{besser}, sections 4 and 5)
for the precise
definitions as well as a full discussion of Coleman integration.

There is a crystalline interpretation of the De Rham fundamental
group via overconvergent connections. That is, we consider the
category $\Un_{an}$ of unipotent overconvergent (iso-)crystals on $Y$.
In fact, there is an equivalence of categories $E:\Un\simeq
\Un_{an}$ (\cite{CL}, Prop. 2.4.1). If we let $y\in Y$, there is the fiber
functor $s_y:\Un_{an}\ra Vect_F$ which associates to a crystal $M$
the horizontal sections on the tube $]y[$ (\cite{besser}, p.26).
Suppose $x\in X$
reduces to $y\in Y$. Then the map $ev_x$ evaluating horizontal
sections at the point $x$ defines an isomorphism of functors
$ev_x:s_y \circ E \simeq e_x $. We can define $s_y$ also for
$y=0,1$ or $\infty$. For example, if $y=0$, then $s_0$ associates
to the overconvergent crystal $M$, a full set of solutions with
coefficients in $\Anlog^a (0)$. A similar discussion holds near $1
$ or $\infty$. For any points $y,y'\in \P^1 (k)$, we have the
group of isomorphisms from $s_{y}$ to $s_{y'}$ giving rise to the
crystalline fundamental groups $\picr (Y,y)$, which are
pro-unipotent algebraic groups over $F$, and the torsors of paths
$\picr (Y, y,y')$ (\cite{besser}, section 3).
Now if $x,x'\in X(F)$ reduce to $y,y'\in Y$,
then we have a natural isomorphism $\pidr(X,x,x')\simeq \picr
(Y,y,y')$ defined by the evaluation isomorphisms $ev_x$ and
$ev_{x'}$. On the other hand, the crystalline torsors have the
natural action of Frobenius endomorphisms $\phi: \picr (Y,y,y')
\ra \picr (Y,y,y')$ (loc. cit.).
That is to say, if $f:\P_k^1 \ra \P_k^1$ is the
geometric Frobenius of $\P_k^1$, then there is a power $f^a$ of $f$
that fixes $y$ and $y'$.  Therefore, the pull-back $(f^a)^*$ defines an
automorphism of the category $\Un_{an}$ that commutes with the
functors $s_y$ and $s_{y'}$. The endomorphism $\picr (Y,y,y') \ra
\picr (Y,y,y')$ induced by any such pull-back will be called a
Frobenius endomorphism and denoted by the same letter $\phi$ when
no danger of confusion is present. Besser and Furusho define fiber
functors $s_w$ on unipotent crystals also for tangential
basepoints $w\in T_y(\P^1)$, $ y\in \P^1\setminus Y$ (\cite{BF}, section 2).
This
construction {\em uses} the equivalence $E: \Un \simeq \Un_{an}$
discussed above. That is, given a unipotent crystal $M$, $\Res_y(M):=E(\Res_x
(E^{-1}(M)))$ where $x=0,1$ or $\infty$ following the value of $y$.
Here, $\Res_x(E^{-1}(M))$ is
 an overconvergent connection on $T_x(\P_F^1)_{an}\setminus
\{0\}$, the analytic space associated to $T_x(\P_F^1)\setminus
\{0\}$. Hence, $\Res_y(M)$ is an overconvergent unipotent crystal on
 $T_y(\P_k^1)\setminus\{0\}$.
Then
$s_w$ associates to $M$, the horizontal sections to  $E(\Res_y
(E^{-1}(M)))$ on the tube $ ]w[$ of $w$ in
$T_x(\P_F^1)_{an}\setminus \{0\}$. If an $F$-rational $v $ in
$T_x(\P_F^1)\setminus \{0\}$ reduces to $w$, then the construction
makes it clear that  evaluation of horizontal sections at $v$
defines a isomorphism of functors $ev_v: s_w\circ E \simeq e_v$.

According to Besser (\cite{besser}, Cor. 3.2) and
Besser-Furusho (\cite{BF}, Thm. 2.8), for any two points $y,y'$,
possibly tangential, we have a unique Frobenius invariant path
$\g_F(y,y')\in \picr(Y,y,y')$. Thus, if $x$ and $x'$ reduce to $y$
and $y'$, we get a Frobenius invariant path in $\pidr(X,x,x')$
which we will denote by $\g_F(x,x')$. Note that if $x$ and $x'$ in
$X(F)$ reduce to the same point $y\in Y$, then this path is
described on connections $(V, \nabla)$ as follows: Given an
element $t_x \in V_x$, let $s$ be the unique horizontal section on
the tube of $y$ such that $s(x)=t_x$. Then
$\g_F(x,x')(t_x)=s(x')$. Now let $v=d/dz$, a tangent vector in
$\P^1$. Then for a point $x \in X(F)$ reducing to $y\in Y$, the
description of the path $\g_F(v,x)\in \pidr(X,v,x)=P(x)$ is
slightly more intricate (\cite{BF}, Prop. 2.11). As above, let $(\bar{V},
\bar{\nabla})$ be the canonical extension of $(V,\nabla)$. Then we
have a canonical isomorphism $\Res_0(V,\nabla)_v \simeq
\bar{V}_0$. An element $t$of $s_0(E(V,\nabla))$ can be identified
with an element of $\Anlog^a (0)\otimes \bar{V}_0$, and hence, can
be written $t=t_0+t_1(\log_a(z))+t_2(\log_a(z))^2+\cdots$ where
the $t_i$ are Laurent series with values in $\bar{V}_0$ converging
in some annulus. The constant term $c_0(t)\in \bar{V}_0$ is then
defined to be the constant term of the Laurent series $t_0$. Now
if we start with a vector $t_v\in \Res_0(V,\nabla)=\bar{V}_0$, we
can find a unique $t\in s_0(E(v,\nabla))$ such that $c_0(t)=t_v$,
and then $\g_F(v,x)(t_v)=\g_F(0,y)(t)(x)$.

 As described in
Deligne (\cite{deligne}, section 12), there is also a canonical element $\g_{DR}(x) \in P(x)$
corresponding to the `De Rham trivialization' of unipotent
connections. So for any $x \in ]Y[$, we get an element $g_x \in U$
such that $ \g_{DR}(x)g_x=\g_F(x)$. We call this map $U_v\Alb:
x\in ]Y[ \mapsto g_x\in U$ the unipotent Albanese map. The
corresponding images $U_v\Alb_n(x) \in U_n$ are p-adic analogues
of the `higher' Albanese maps of Hain \cite{hain}, except for the fact that in
the p-adic case, we can dispense of periods using the
Frobenius invariant path. We note that if we choose a
different base-point $z$, then there are obviously Albanese maps
$U_z\Alb$ and $U_z\Alb_n$ defined analogously to $U_v\Alb$.

 Let
$V:=F<<A,B>>$ be the ring of non-commutative power series over $F$
in the variables $A$ and $B$. $A$ and $B$ act as endomorphisms of
the fiber functor of evaluation at $v$ via the residue of a
connection. $U$ can be identified with the group-like elements of
$V$ while $U_n$ is identified with the group-like elements of
$V_n:=V/I^{n+1}$ (\cite{deligne}, 16.1.4 and Prop. 16.4).

 Now consider the trivial pro-vector bundle $\cV$ on $X$ with
fiber $V$ and connection
$$\nabla (1)=Adz/z+Bdz/(z-1)$$

As in Furusho (\cite{furusho}, Thm. 3.4), there is a unique horizontal section $G_a(z)$ with
coefficients in the Coleman functions $\Col_a$ such that
$G_a(z)\ra z^{A}:=\exp (\log_a(z)A) $ as $z\ra 0$. This means that
$G_a(z)z^{-A}$ is of the form $1+u(z)$ where $u(z) \in \An
(]0[)\otimes V$ and $u(0)=0$. Write $G_a(z)=\Sigma_wg^a_w(z)w$
where $w$ runs over the words in $A,B$. In fact, $G_a(z) \in U$
(i.e., it is group-like) for each $z$ (\cite{furusho}, proof of
Prop. 3.39).
\begin{prop} For any $z\in ]Y[$,
$UAlb (z)=G_a(z):=r_y(G_a)(z)$ where $y$ is the reduction of $z$.
\end{prop}

{\em Proof.} The representation of $U$ associated to the bundle
$\cV$ is simply $V$ with the canonical action of $U$ by left
multiplication. The De Rham trivialization assigns to any point
$z$, the isomorphism $V_v\simeq V\simeq V_z$. We know that
$G_a(z)$ is group-like for each $z$. If  we expand $G_a(z)$ as a
power series in $\log_a(z)$, its constant term is 1. Hence, given
any vector $l\in V$, $G_a(z)l$ is a horizontal section of $\cV$
such that its constant term near zero is $l$. Furthermore, it is a
horizontal section of $\cV$ with coefficients in the Coleman
functions. As described above, according to Besser and Furusho,
taking the constant term of a horizontal section near zero
corresponds to the fiber functor associated to the tangential
base-point $v=d/dz$. It is essentially tautological then that
$g_F(z)l=G_a(z)l$, but let us briefly sketch the logic: To
construct $G_a(z)$, one constructs the solution  $t=\Sigma t^w w$
to the K-Z equation satisfying the asymptotic condition in
$\An(]0[)[\log_a(z)] \otimes V$ and then translates it to a
solution $t_y=\Sigma t^w_y w$ in $\Anlog^a(y)\otimes V$ using the Frobenius
invariant isomorphism $\g_F(y): (\Anlog^a(0)\otimes V)^{\nabla=0}
\simeq (\Anlog^a (y)\otimes V)^{\nabla=0}$. One gets thereby a
Coleman function $g^a_w=[(\cV, f_w, \{t_y\})]$ for each word $w$,
where $f_w$ denote the projection to $\O^{\dagger}_X$ given by the
$w$-component. One then gets the solution $[G_a]:=\Sigma g^a_w w$
which has values in $V$ and coefficients in $Col_a$. From the
construction, we have $r_y(g^a_w)=t_y^w$ for each $w$ and $y$.
Now, if $l\in V$, by the construction above,  we get that
$\g_F(z)l=t_y(z)l=r_y(G_a)(z)l=G_a(z)l$. Since the action of $U$
on $V$ is faithful, we get $\g_F(z)=G_a(z)$. $\Box$

\begin{thm}

The functions $g^a_w(z)$ on $]Y[$ are linearly independent over
$F$.
\end{thm}
{\em Proof.} Given two different choices of branches $a, b\in F$
for the p-adic logarithm, there is an isomorphism of rings $Col_a
\simeq Col_b$ characterized by $f(z)\mapsto f(z)$ for $f(z)\in
A^!$ and $\log_a(z) \mapsto \log_b(z)$. The proof will use the
comparison between $G_a(z)$ and $G_b(z)$. First, note that for
functions in $\Col_a$, linear independence over $F$ is equivalent
to independence over the ring of locally constant functions. To
see this,  suppose $f_i$ are Coleman functions linearly
independent over the constants and assume that $\Sigma_ic_i f_i=0$
where the $c_i$ are locally constant. For any point $z$ in $]Y[$,
there exists a neighborhood $z\in O$ such that the $c_i$ are
constant on $O$. Then $\Sigma_i c_i(z)f_i=0$ on $O$. By the
identity principle for Coleman functions, this implies
$\Sigma_ic_i(z)f_i=0$ everywhere. Then, using the linear
independence over the constants, this implies that $c_i(z)=0$ for
all $i$.

Using the differential equation, it is  easy to see that
$G_a=G_bc(z)$ for a locally constant function $c(z)$. By the
asymptotic condition at 0, we get that $c(z)\approx z_b^{-A}z_a^A$
as $z\ra 0$, where we write the subscripts $a,b$ to denote the
dependence on the different logs. In fact, let's write this out a
bit more precisely.

We have that $G_a(z)z_a^{-A} \approx 1$ and $G_b(z)z_b^{-A}
\approx 1$, so from, $G_a(z)z_a^{-A}=G_bz_b^{-A}z_b^Ac(z)z_a^{-A}$
we get
$$z_b^Ac(z)z_a^{-A}\approx 1$$
or more precisely,
$$z_b^Ac(z)z_a^{-A}=1+u_A(z)A+u_B(z)B+\cdots$$
where $u_A(z)$ etc. are rigid analytic functions that vanish at
the origin. So
$$c(z)=1+(\log_a(z)-\log_b(z)+u_A(z))A+u_B(z)B+\cdots$$
By local constancy, we get $u_A=0, u_B=0$. That is,
$c(z)=1+(\log_a(z)-\log_b(z))A+d(z)$ where $d(z)$ involves only
words of length $\geq 2$.

Recall that there are also functions $G_a^1$ and $G_b^1$
satisfying the same differential equations and the asymptotic
condition
$$G_a^1\approx (1-z)_a^B, G_b^1\approx (1-z)_b^B$$near
as $z\ra 1$. By an argument entirely similar to that above, this
gives us $G_a^1(z)=G_b^1(z)c^1(z)$ for a locally constant function
$$c^1(z)=1+(\log_a(1-z)-\log_b(1-z))B +d^1(z)$$ where $d^1(z)$
also involves only words of length $\geq2$. Meanwhile, the p-adic
Drinfeld associator $\Phi$relates the two asymptotics, so that
$G_a(z)=G^1_a(z)\Phi$ and $G_b(z)=G^1_b(z)\Phi$ for the same $\Phi
\in C$. This gives us the additional relation
$$G_a(z)=G^1_a(z)\Phi=G_b^1(z)c^1(z)\Phi=G_b(z) \Phi^{-1}c^1(z)\Phi$$
 or
$c(z)=\Phi^{-1}c^1(z)\Phi$ near $z=1$. The explicit formula of
Furusho (example 3.35) shows that $\Phi$ involves no linear terms,
so that $c(z)=1+(\log_a(1-z)-\log_b(1-z))B+d'(z)$ near $z=1$,
again for $d'$ with words of length $\geq 2$. Write
$$G_a(z)=\Sigma g_w^a(z)w, \ \ G_b(z)=\Sigma g_w^b(z)w$$
Now we will prove that the $g^a_w$ are linearly independent by
induction on the length $n=l(w)$ of $w$. The statement is clearly
true for $n=1$. Assume that the $g^a_w$ are linearly independent
for $w$ of length $\leq n-1$.

Suppose $\Sigma a_wg_w^a=0$ for some $a_w$ running over $w$ of
length $\leq n$. Then $\Sigma_wa_wg_w^b=0$, that is, linear
relations are preserved under a change of logs. Write the
relations as
$$\begin{array}{rl}\Sigma_{l(w)=n-1}a_{wA}g^a_{wA}+&
\Sigma_{l(w)=n-1}a_{wB}g^a_{wB}\\
+\Sigma_{l(w)= n-2}
a_{wA}g^a_{wA}(z)+&\Sigma_{l(w)=n-2}a_{wB}g^a_{wB}\\
+\Sigma_{l(w)\leq n-2}a_{w}g^a_{w}(z)=0\end{array}$$ and
$$\begin{array}{rl}\Sigma_{l(w)=n-1}a_{wA}g^b_{wA}+&\Sigma_{l(w)=n-1}a_{wB}g^b_{wB}\\
+\Sigma_{l(w)= n-2}
a_{wA}g^b_{wA}(z)&+\Sigma_{l(w)=n-2}a_{wB}g^b_{wB}\\
+\Sigma_{l(w)\leq n-2}a_{w}g^b_{w}(z)=0 \end{array}$$ Using the
relation between $G_a$ and $G_b$ near $z=0$, we get from the first
equation
$$\begin{array}{rl}
\Sigma_{l(w)=n-1}(a_{wA}(g^b_{wA}(z)+(\log_a(z)-\log_b(z))g^b_{w}(z))
+&\Sigma_{l(w)=n-1}a_{wB}g^b_{wB} \\ +\Sigma_{l(w)= n-2}
a_{wA}(g^b_{wA}+(\log_a(z)-\log_b(z))g^b_w(z))
+&\Sigma_{l(w)=n-2}a_{wB}g^b_{wB}\\+\Sigma_{l(w)\leq
n-2}e_{w}g^b_{w}(z)=0
\end{array}$$
and subtracting this from the second equation, we get

$$(\log_a(z)-\log_b(z))\Sigma_{l(w)=n-1}a_{wA}g^b_{w}(z)
+\Sigma_{l(w)\leq  n-2} e'_wg^b_w(z)=0$$ From which we get
$a_{wA}=0$ when $l(w)=n-1$. Similarly, $a_{wB}=0$ for $l(w)=n-1$.
Then by induction, we get $a_w=0$ for all $w$. $\Box$

We will now suppress the choice of log from the notation. Suppose
we take another basepoint $x$ instead of $v$. We wish to compare
$U_x\Alb$ and $U_v\Alb$. Using the uniqueness of the Frobenius
invariant path we see that $U_x\Alb (y)= (\g_{DR}(x))U_v\Alb
(y)\g_F^{-1}(x)$. If $\alpha_w$ is the function on $U$ that picks off
the coefficient of the word $w$, then the formula for the
comultiplication on $U$ shows that $$\alpha_w(U_v\Alb
(\cdot)\g_F^{-1}(x))=\alpha_w(U_v\Alb
(\cdot))+\Sigma_{\{w'| w=w'w''\}}c^w_{w'}\alpha_{w'}(U_v\Alb (\cdot))$$
where the $c^w_{w'}$ are constants (depending on $x$).
Thus, by induction on the length of $w$, the
linear independence of the $g_w$ implies the linear independence
of the $\alpha_w(U_v\Alb (\cdot)\g_F^{-1}(x))$. Similarly, by
considering multiplication on the other side,
$g_w(x,\cdot):=\alpha_w((\g_{DR}(x))U_v\Alb (\cdot)\g_F^{-1}(x)$ are
linearly independent.

Although we will not need it, as discussed in Furusho \cite{furusho}, the
function $G_a(z)$ extends to a function on $X(\C_p)$
and the same argument as that given above shows that the coefficients
are linearly independent over $\C_p$.

\section{The non-abelian method of Chabauty}

Now let $\cX=\P^1_{\Z [1/S]}-\{0,1,\infty\}$. Let $X$ be the
generic fiber of $\cX$ and $Y$ the special fiber over a prime
$p\notin S$. Let $\bar{x}$ be a $\Z_p$ point of $\cX$ and $x\in X$
its restriction to $X\otimes \Q_p$. $y\in Y$ denotes its reduction
to $\F_p$.
 Associated to $X\otimes \Q_p$
we have the $\Q_p$-unipotent \'etale fundamental group
$\piet(X\otimes \Q_p,x)$ and the De Rham fundamental group $\pidr
(X\otimes Q_p,x)$ which is a unipotent group over $\Q_p$. We
denote the corresponding coordinate rings by $\Aet$ and $A_{DR}$.
Let $G=G_{\Q_p}$ be the Galois group of $\Q_p$. We also denote by
$[\Aet]_n$ and $[A_{DR}]_n$ the subalgebras corresponding to the
lower central series of the fundamental groups, where we
 set the index so that the abelianization corresponds to $n=1$.
In the following, if a definition, statement, or proof
of a statement works the same way for a group or
a quotient in the descending central series,
we will give it only for one of them unless a danger of
confusion presents itself. In fact, the only cause for concern
arises from  the p-adic comparison
theorem.
According to Vologodsky (\cite{vologodsky}, 1.9, Thm A), $[\Aet]_n$
is a De Rham representation of $G$ for $(p-1)/2\geq n+1$ and
$[A_{DR}]_n\simeq D([\Aet]_n):=([\Aet]_n)\otimes_{\Q_p}B_{DR})^G$
together with the induced weight and Hodge filtration. Also,
$[\Aet]_n\otimes B_{DR}\simeq [A_{DR}]_n \otimes B_{DR}$ compatibly with
the all the structures. Here of course, the Hodge filtration on
$\Aet$ is trivial ($F^0\Aet=\Aet, F^1\Aet=0$), the weight
filtration on $B_{DR}$ is trivial ($W_0B_{DR}=B_{DR},
W_{-1}B_{DR}=0$), and the Galois action on $A_{DR}$ is trivial. For
$\Q_p$-algebras $R$ of finite type, we
will be considering  torsors $P$ of $[\piet^R]_n:=[ \piet(X,x)]_n\otimes_{\Q_p}R$
of {\em De Rham type}. (We will also simply call these `De Rham torsors.)
We first define it for finite field extensions $L$ of $\Q_p$
by the condition that the
coordinate ring $\cP$ of the $\piet^L$-torsor $P$
has the property that the natural
inclusion $D(\cP) \ra \cP\otimes B_{DR}$ induces an isomorphism
$D(\cP)\otimes B_{DR}\simeq \cP\otimes B_{DR} $, where
$D(\cP):=(\cP\otimes B_{DR})^{G}$. Here, we assume that the
isomorphism respects all the structures, namely, the Galois
action, the weight filtration, and the Hodge filtration. Since
each level $W_n$  of the weight filtration is finite-dimensional,
the De Rham property is equivalent to the condition $\dim
W_n(D(\cP))=\dim W_n\cP$ for each $n$.
Now for $R$ of finite type, say that $P$ is De Rham if its
base change to all the closed points of $\Spec (R)$ is De Rham.
Denote by $\piet^B$, the
base change $\piet^B:=\Spec(\Aet\otimes_{\Q_p} B_{DR})$. The constructions
of section 1 allow us to put a canonical structure of an affine
variety on the cohomology set $H^1(G, [\piet]_n)(\Q_p)$ (for $B=F=K=\Q_p$
with trivial action). The map 
 $H^1(G,
[\piet]_n) \ra H^1 (G, [\piet^B]_n)$   also defines a
subset $H^1_{f}(G, [\piet]_n)(\Q_p)\subset H^1(G, [\piet]_n)(\Q_p)$
as the inverse image of the basepoint, which in fact is the
set of $\Q_p$ points of a subvariety $H^1_{f}(G, [\piet]_n)
\subset H^1(G, [\piet]_n)$ as discussed in section 1.

\begin{prop} Assume $(p-1)/2\geq n+1$.
 Let $L$ be a finite extension of $\Q_p$ and
let $$[P]\in H^1(G, [\piet]_n)(L)$$ be the cohomology class of the
torsor $P$. Then $P$ is De Rham if and only  if $[P]\in H^1_{f}(G,
[\piet]_n)(L)$.
\end{prop}
{\em Proof.} Suppose the class $[P]$ becomes trivial in $H^1 (G,
[\piet^B]_n)(L)$. This means $\cP\otimes B_{DR} \simeq \Aet \otimes_{\Q_p}
 L \otimes_{\Q_p}
B_{DR}$ as Galois representations with weight filtrations. But
then, $D(\cP)\simeq D([\Aet]_n \otimes L)=D([\Aet]_n)\otimes L$ so
that $\dim_L W_mD(\cP)=\dim_{L} W_mD([\Aet]_n)\otimes L=\dim_{L}
W_m[\Aet]_n\otimes L=\dim_L\dim W_m \cP$. So $P$ is De Rham.
Conversely, assume $\cP\otimes B_{DR}\simeq D(\cP)\otimes B_{DR}.$
$D(\cP)$ is a torsor for $D([\Aet]_n\otimes L)$. By unipotence,
there is a point $x\in D(P)$ inducing an isomorphism $$D(\cP)
\simeq D([\Aet]_n\otimes L)$$ Also, we have an equality of
dimensions
$$\dim W_mD(\cP)=\dim W_m \cP=\dim W_m([\Aet]_n\otimes L)=\dim
(W_m D([\Aet]_n))\otimes L$$ Thus, $D([\Aet]_n\otimes L)\simeq
D(\cP)$ together with the weight filtration. Since the Galois
actions are trivial on $D([\Aet]_n \otimes L)$ and $D(\cP)$, we
get
$$\cP\otimes B_{DR}\simeq D(\cP)\otimes B_{DR} \simeq D([\Aet]_n\otimes L)
\otimes B_{DR} \simeq
[\Aet]_n\otimes L \otimes B_{DR}$$ as torsors with Galois action. Therefore,
$[\cP\otimes B_{DR}]$ is trivial. $\Box$

 Recall from the previous section that we have a Frobenius map $$\phi:\pidr=\pidr (X\otimes \Q_p,x) \ra
 \pidr(X\otimes \Q_p,x)
 $$induced from the comparison with the crystalline fundamental
 group.

If $R$ is a finitely-generated $\Q_p$-algebra, by a torsor for
$\pidr\otimes R$, we will mean an R-scheme $P=\Spec(\cP)$ equipped
with a non-negatively indexed decreasing (Hodge) filtration
$F^i\cP$, an increasing weight filtration $W_n\cP$ by finitely
generated $R$-submodules, and an $R$-algebra (Frobenius)
automorphism $\phi_P: \cP \ra \cP$. $P$ is equipped with an action
of $\pidr\otimes R$ that gives it the structure of a torsor for
the pro-unipotent group $\pidr$ in the usual sense. However, we
require that the structure map $a: \cP \ra \cP \otimes A_{DR}$ for
the action preserves both filtrations and the Frobenius map. For
any maximal ideal $m$ of $R$, we get a torsor $P\otimes R/m$ for
$\pidr\otimes R/m$. We will say the family is admissible if for
each of these fibers, $P\otimes R/m$ is of the form $D(P_{et}(m))$
for a De Rham $\piet^{R/m}$-torsor $P_{et}(m)=\Spec
(\cP_{et}(m))$. Note that we have $\dim W_n \cP\otimes R/m=\dim
W_n \cP_{et}(m) =\dim W_n \Aet\otimes R/m$ is constant over
maximal ideals $m$, and hence, $W_n\cP$ is locally-free of
constant rank equal to this dimension. We can extend the
definition above in an obvious way to torsors for
$[\pidr]_n\otimes R$.

According to Besser (\cite{besser}, Thm 3.1),
the Lang map $L:\pidr \ra \pidr$ that sends
 $g$ to
$g^{-1}\phi(g)$ is an isomorphism. Given any finitely-generated
$\Q_p$-algebra $R$ and a torsor $P_{DR}$ for $\pidr\otimes R$ this
implies (\cite{besser},  proof of Cor. 3.2)that there exists a
unique Frobenius invariant element $p_F(R)$ in $P_{DR} (R)$. The
filtration $F$ on the coordinate ring $\cP_{DR}$ determines a
filtration by subschemes $F^iP_{DR}$. $F^0P_{DR}$ in particular is
the subscheme defined by the ideal $F^1\cP_{DR}$.

\begin{lem} Let $P_{DR}$ be a torsor for
$\pidr \otimes R$  (or $[\pidr]_n\otimes R$). If $P_{DR}$ is
admissible, then there is a unique element $p_{DR}(R)$ in
$F^0P_{DR}(R)$.
\end{lem}
{\em Proof.} If $R$ is a field, this follows from Deligne (\cite{deligne},
Prop. 7.12) and the
fact that both the filtrations F on $\cP_{DR}$ and on
$A_{DR}$ are gradable (by
Wintenberger \cite{wintenberger}). That is, this implies that $P_{DR}$ is associated
to  $F^0P_{DR}$ which is a
$F^0\pidr$-torsor, and $F^0\pidr$ is the trivial
group (\cite{deligne}, section 12),
so $F^0P_{ DR}$ is scheme-theoretically a point.
 Now,
$P_{DR}\otimes R/m$ is a torsor for $\pidr\otimes R/m$ for any
maximal ideal $m$ of $R$. Consider the structure map $R\ra
\cP_{DR}/F^1\cP_{DR}$. By the result over fields, we see that $R/m
\ra \cP_{DR}/F^1\cP_{DR}\otimes R/m $ is an isomorphism for each
maximal ideal $m$. Hence, $R\simeq \cP_{DR}/F^1\cP_{DR}$, the
inverse of which provides the element of $F^0P_{DR}$. Uniqueness
is obvious from this construction. $\Box$

Define the variety $Z^1_f(G, [\piet]_n)$ to be the inverse image of
$H^1_{f}(G, [\piet]_n)$ under the natural projection $Z^1(G, [\piet]_n)
\ra H^1(G, [\piet]_n)$. If we consider the map $$Z^1(G, [\piet]_n)(R) \ra
Z^1(G, [\piet^B]_n)(R),$$ $Z^1_f(G, [\piet]_n)(R)$ consists exactly of
those cocycles whose image in $Z^1(G, [\piet^B]_n)(R)$ lie in the
image of the map $$[\piet^B]_n(R) \ra Z^1(G, [\piet^B]_n)(R).$$
 Taking the
fiber product, we get a surjective map of functors
$$[\cZ]_n:= [\piet^B]_n \times_{Z^1(G, [\piet^B]_n)} Z^1_f(G, [\piet]_n)
\ra Z^1_f(G, [\piet]_n)$$ Now, since the map $[\piet^B]_n(R) \ra
Z^1(G, [\piet^B]_n)(R)$ takes $u$ to $ug(u^{-1})$, we easily see
that $[\cZ]_n$ has the structure of torsor for the functor
$[\piet^B]_n^G$ whose value on $R$ is just the $G$-fixed part of
$[\piet(B\otimes R)]_n$. That is we need to consider the
$G$-invariants in the $B$-algebra homomorphisms from
$[\Aet]_n\otimes B $ to $ B\otimes R$. Now, let $(p-1)/2\geq n+1$.
Then $[\Aet]_n\otimes B= [A_{DR}]_n \otimes B$ and $[A_{DR}]_n$
has trivial $G$-action. So
$$\AlgHom_B ([\Aet]_n\otimes B, B\otimes R)=\AlgHom_B ([A_{DR}]_n \otimes B, B\otimes R)
=\AlgHom_{\Q_p}([A_{DR}]_n, B\otimes R)$$ Hence, the $G$-invariants
are just
$$\AlgHom_{\Q_p} ([A_{DR}]_n, (B\otimes R)^{G})=\AlgHom_F ([A_{DR}]_n, R)$$
That is, $[\cZ]_n$ is a torsor for $[\pidr]_n$ over $Z^1_f(G,
[\piet]_n)$. Now, pulling back to $H^1_f(G, [\piet]_n)$ via a
section $H^1_f(G, [\piet]_n) \ra Z^1_f(G, [\piet]_n)$, we get a
$[\pidr]_n$-torsor over $H^1_f(G, [\piet]_n)$, which we will also
denote by $[\cZ]_n$. This torsor has two sections $\g_{DR}\in
F^0[\cZ]_n (R)$ and $\g_{F}\in [\cZ(R)]_n^{\phi=1}$ where $R$ is
now the coordinate ring of $H^1_f(G, [\piet]_n)$. The transporter
in $[\pidr]_n (R)$ between these two points gives us an $R$-point
of $[\pidr]_n$, and hence, an algebraic map $\cD: H^1_f(G_F,
[\piet]_n) \ra [\pidr]_n$. Given a class $c\in H^1_f(G_F,
[\piet]_n)$, there is the torsor $P_c=\Spec(\cP(c))$ corresponding
to it and a corresponding $[\pidr]_n$-torsor $D(P_c)$. By
considering the transporter between $F^0D(P_c)$ and
$D(P_c)^{\phi=1}$, we get an element $c' \in [\pidr]_n$.

\begin{lem} $c'=\cD (c)$.
\end{lem}
{\em Proof.} It suffices to check that $D(P_c)$ coincides with the
torsor $([\cZ]_n)_c$. We check this on points (for arbitrary $R$).
$([\cZ]_n)_c(R)$ consists of $u\in \piet(B\otimes R)$ such that
$ug(u^{-1})=c(g)$ for all $g\in G$. This is equivalent to
$u=c(g)g(u)$. Now we examine the points of $D(P_c)$. These are the
$G$-ivariant homomorphisms $x:\cP(c) \ra B\otimes R$. Since
$\cP(c)$ is just $[\Aet]_n\otimes R$ with  the action twisted by
$c$, this is equivalent to looking at the $G$-invariant
homomorphisms $x:[\Aet]_n\otimes B\otimes R \ra B\otimes R$. Such
homomorphisms are in particular points in $[\piet]_n(B\otimes R)$.
Let us impose the invariance condition. For this, recall also the
formula $c(g^{-1})=g^{-1}c(g)^{-1}$ for a 1-cocycle $c$. Now the
$c$-twisted $G$-action takes $x$ to the homomorphisms whose value
on $a$ is
$$g(x(\phi_{c(g^{-1})}(g^{-1}(a))))=g((g^{-1}(a))(c(g^{-1})^{-1}x))=
g(g^{-1}(a(g(c(g^{-1})^{-1})g(x))))=a(c(g)g(x))$$
that is, the point $c(g)g(x)$. Thus, $G$-invariance is the
same as $c(g)g(x)=x$, which is exactly the condition for
$([\cZ]_n)_c(R)$.
$\Box$

{\em Proof of Siegel's theorem.} Let $\G:=\Gal(\bar{\Q}/\Q)$ and
consider a ring $\Z[1/S]$ of $S$-integers. Choose a $p\notin S$
and let $T=S\cup \{p\}$. Below, we will put further restrictions
on $p$. We will consider the pro-unipotent completion of the
$p$-adic \'etale fundamental group of $\bar{X}=(\P^1\setminus
\{0,1,\infty \})\otimes \bar{\Q}$ and denote it by $\piet$. Here
we choose an $S-$integral point $x$ as a base-point, if it exists
(otherwise, we are done). Thus, we get an action of $\G_T$, the
Galois group of the maximal extension of $\Q$ unramified outside
$T$, on $\piet$. Let $G$ be a decomposition group at $p$. We can
then restrict the representation to $G$, which can also be
interpreted as the action on the $\Q_p$-unipotent fundamental
group of $(\P^1\setminus \{0,1,\infty \})\otimes \bar{\Q}_p$.

Any other integral point $y$ determines a point $[[\piet (\bar{X},
x,y)]_n]$ in $H^1(\G_T,[\piet]_n)(\Q_p)$, which we regard as an
affine variety following the construction of section 1
(again with $K=F=B=\Q_p$). In
fact, if we define $H^1_f(\G_T, [\piet]_n)$ to be the inverse
image of $H^1_f(G,[\piet]_n)$ with respect to the restriction map
$$H^1(\G_T,[\piet]_n)\ra H^1(G,[\piet]_n)$$
the class of the point $y$ lies in $H^1_f(\G_T,\piet)(\Q_p)$ for
$(p-1)/2\geq n+1$ (\cite{vologodsky}, 1.9, Thm. A). There is a
commutative diagram
$$\begin{array}{ccc}
\cX(\Z[1/S])& \ra &X(\Q_p)\cap ]Y[\\
\da & & \da \\
H^1_f(\G_T,\piet)(\Q_p) & \ra & H^1_f(G,\piet)(\Q_p)
\end{array}
$$
simply from the compatibility with base change of the fundamental
groups and torsors of paths while loc. cit. yields the commutative
diagram
$$\begin{array}{ccc}
X(\Q_p) \cap ]Y[& \ra & [\pidr]_n (\Q_p)\\
\da & & \da \scriptstyle{=}\\
H^1_f(G,[\piet]_n)(\Q_p) & \ra & [\pidr]_n(\Q_p)
\end{array}
$$
for $(p-1)/2 \geq r+1$, where the upper horizontal arrow is of
course the level n Albanese map $UAlb_n$. Hence, $UAlb_n (X(\Q_p)
\cap ]Y[)$ lies inside the image of $H^1_f(\G_T,[\piet]_n)(\Q_p)$.

But there is a bound on the dimension of the latter. That is, we
have the exact sequences
$$0\ra H^1_f(\G_T,\piet^{n}/\piet^{n+1})(\Q_p) \ra H^1_f(\G_T,[\piet]_n)(\Q_p)\ra
H^1_f(\G_T,[\piet]_{n-1})(\Q_p)\ra 0$$ in the sense that the
middle term is a fibration with fibers isomorphic to the first
term, and an analogous one for the De Rham fundamental group:
$$0\ra  \pi_{DR}^{n}/\pidr^{n+1}\ra [\pidr]_{n+1}\ra [\pidr]_n\ra 0$$
(\cite{deligne}, Prop. 16.3, 16.4) On the other hand, we have
$\piet^{n}/\piet^{n+1}\simeq \Q_p(n)^{r_n}$ and
$\pi_{DR}^{n}/\pi_{DR}^{n+1}\simeq \Q_p(n)^{r_n}$ for some rank
$r_n$ growing to infinity and, importantly, independent of the
choice of $p$. So
$$\dim \pi_{DR,n}=r_1+r_2+\cdots r_n.$$ But by Soul\'e's vanishing
theorem (\cite{soule}), $H^1(G_T, \Q_p (2n))=0$ for each positive $n$, while
$H^1(G_T, \Q_p(2n+1))$ has dimension 1 for $n\geq 3$. Hence,
$$\dim H^1_f(\G_T,(\piet)_n)(\Q_p) \leq R+r_3+r_5+\cdots+r_{2\lfloor n/2 \rfloor}$$
where $R=2\mbox{rank}(\Z[1/S]^*)$. Therefore, for $n$ sufficiently
large, $\dim H^1_f(\G_T,[\piet]_n) < \dim [\pidr]_n$. Now we
choose $p$ such that $(p-1)/2\geq n+1$ so that we have the
commutative diagrams above, and we get that some element of the
coordinate ring of $[\pidr]_n$ must vanish on the image of $\dim
H^1_f(\G_T,[\piet]_n)(\Q_p) $. But when pulled back to $X(\Q_p)
\cap ]Y[$ via the unipotent Albanese map, the elements of this
coordinate rings are exactly linear combinations of the
$g_w(x,\cdot)$ of section 2. Thus, there is a non-zero Coleman
function that vanishes on $\cX(\Z[1/S])$. Since $X(\Q_p)\cap ]Y[$
is compact, this implies the desired finiteness. $\Box$.

\medskip
{\bf Acknowledgements:} I seem to have accumulated quite a debt of
gratitude even in writing this simple paper: (1) to Dick Hain for
teaching me about iterated integrals; (2) to Jean-Marc Fontaine
and Luc Illusie for many conversations and lectures on p-adic
Hodge theory as well as an invitation to Orsay in September 2003;
(3) to Hidekazu Furusho for patiently explaining to me his results
on $p$-adic multiple polylogs; (4) to Makoto Matsumoto, who
generously supported my visit to Kyoto and Hiroshima from October
to December, 2003; (5) to Akio Tamagawa, who also arranged the
visit and listened patiently to a premature version of the proof
presented here; (6) to Ofer Gabber who helped me clarify the proof
of Prop. 2; (7) to Gerd Faltings for his interest, encouragement,
 and helpful suggestions;
(8) and to the unbelievably hospitable staff in
the international office at RIMS, especially Ms. Eriko Edo. But my
deepest gratitude goes to Shinichi Mochizuki, who explained to me
his profound ideas relating fundamental groups to Diophantine
problems in innumerable hours of stimulating conversation, and
without whose strong encouragement I would never have thought
through and written down the details of this paper. More
specifically, he corrected  an earlier naive approach to the map `$C$' by
emphasizing to me the importance of non-linearity.

{\footnotesize DEPARTMENT OF MATHEMATICS, UNIVERSITY OF ARIZONA,
TUCSON, AZ 85721, U.S.A. }

{\footnotesize EMAIL: kim@math.arizona.edu}

\end{document}